\documentclass{amsart}
\usepackage{graphicx}
\newtheorem{thm}{Theorem}[section]
\newtheorem{cor}[thm]{Corollary}

\theoremstyle{definition}

\theoremstyle{remark}
\newtheorem{rem}[thm]{Remark}
\numberwithin{equation}{section}

\begin{document}

\title[Left Derivations and Strong Commutativity Preserving Maps]
{Left Derivations and Strong Commutativity Preserving Maps on Semiprime $\Gamma$-Rings}%
\author{Xiaowei Xu}%
\address{Xiaowei Xu: College of Mathematics, Jilin University, Changchun 130012, PR
     China}%
\email{xuxw@jlu.edu.cn}%

\author{Jing Ma}%
\address{Jing Ma: College of Mathematics, Jilin University, Changchun 130012, PR
     China}%
\email{jma@jlu.edu.cn}%

\author{Yuan Zhou}%
\address{Yuan Zhou: College of Mathematics, Jilin University, Changchun 130012, PR
     China}%
\email{zhouyuan150630@126.com}%

\thanks{The paper is supported by the NNSF of China (No. 10871023 and No.11071097),
  211 Project, 985 Project and the Basic Foundation for Science Research from Jilin University.}%

\begin{abstract}
In this paper, firstly as a short note,  we prove that a left
derivation of a semiprime $\Gamma$-ring $M$ must map $M$ into its
center, which improves  a result by Paul and Halder and some results
by Asci and Ceran. Also we prove that a semiprime $\Gamma$-ring with
a strong commutativity preserving derivation on itself must be
commutative and that a strong commutativity preserving endomorphism
on a semiprime $\Gamma$-ring $M$ must have the form
$\sigma(x)=x+\zeta(x)$ where $\zeta$ is a map from $M$ into its
center, which extends some results by Bell and Daif to semiprime
$\Gamma$-rings.
\end{abstract}
\maketitle
\section{Introduction}
In 1964,  Nobusawa \cite{Nobusawa-1964-OsakaJM-GammaNRing} had
introduced the notion of a $\Gamma$-ring, which was extended by
Barnes \cite{Barnes-1966-PJM-GammaRing} in 1966 so that the quotient
$\Gamma$-ring of a $\Gamma$-ring can be defined reasonably. At
present the notion by Nobusawa is called a $\Gamma_N$-ring, and the
notion by Barnes is called a $\Gamma$-ring. A $\Gamma_N$-ring is a
$\Gamma$-ring, and there exists a $\Gamma$-ring being not a
$\Gamma_N$-ring.  Barnes \cite{Barnes-1966-PJM-GammaRing} has
defined a $\Gamma$-ring as following: Let $M$ and $\Gamma$ be two
additive abelian groups. If there exists a map $(a,\alpha,b)\mapsto
a\alpha b$ of $M\times \Gamma\times M\rightarrow \Gamma$ satisfying
the conditions
\begin{itemize}
\item $(a+b)\alpha c=a\alpha c+b\alpha c$, $a(\alpha+\beta)b=
a\alpha b+ a\beta b$, $a\alpha(b+c)= a\alpha b+ a\alpha c$,
\item $(a\alpha b)\beta c=a\alpha(b\beta c)$
\end{itemize}
for all $a,b,c\in M$ and $\alpha,\beta\in \Gamma$, then $M$ is
called a $\Gamma$-ring.

 An associative ring $R$ can be seen as a $\Gamma$-ring. For example
 it is pointed out explicitly in
\cite{Chakraborty-Paul-2010-BIMS}  that an associative ring $R$ is a
$\Gamma$-ring with $\Gamma=U$ or $\mathbb{Z}$ where $\mathbb{Z}$ is
the ring of integers and $U$ is an ideal (even an additive subgroup)
of $R$. Some properties of $\Gamma$-rings as contrast to those of
general rings have also been obtained by Barnes
\cite{Barnes-1966-PJM-GammaRing}, Kyuno
\cite{Kyuno-1978-PJM-PrimeGammaRing} and Luh
\cite{Luh-1969-MichiganMJ-Simple-Gamma}.

 An additive subgroup $U$ of a $\Gamma$-ring $M$ is called a left
(resp. right) ideal of $M$ if $M\Gamma U\subseteq U$ (resp. $U\Gamma
M\subseteq U$). A left ideal $U$ of a $\Gamma$-ring $M$ is   called
an ideal of $M$ if it is also a  right ideal of $M$. The set
$Z(M)=\{a\in M|a\alpha b=b\alpha a, \forall b\in M, \forall
\alpha\in\Gamma\}$ is called the center of $M$.

A $\Gamma$-ring $M$ is called prime if $a\Gamma M\Gamma b=0$ with
$a,b\in M$ implies $a=0$ or $b=0$.  A $\Gamma$-ring $M$ is called
semiprime if $a\Gamma M\Gamma a=0$ with $a\in M$ implies $a=0$. The
notion of a (resp. semi-)prime  $\Gamma$-ring is an extension for
the notion of a (resp. semi-)prime  ring.

Recall that an additive map $\delta$ from a ring $R$ into itself is
called a left derivation if $\delta(xy)=x\delta(y)+y\delta(x)$ holds
for all $x,y\in R$. In 1990, Bre\v{s}ar and Vukman
\cite{Bresar-Vukman-1990-PAMS-left} firstly introduced the notion of
a left derivation in a ring and proved that a left derivation of a
semiprime ring $R$ must map $R$ into its center.

Similarly an additive map $\delta$ from a $\Gamma$-ring $M$ into
itself is called a left  (resp. right) derivation if $\delta(x\alpha
y)=x\alpha\delta(y)+y\alpha\delta(x)$  (resp. $\delta(x\alpha
y)=\delta(y)\alpha x+\delta(x)\alpha y$)  holds for all $x,y\in M$
and $\alpha\in\Gamma$.

Naturally the notions of  derivations and endomorphisms for  rings
have extensions  for  $\Gamma$-rings. An additive map $\mu$ from a
$\Gamma$-ring $M$ into itself is called a derivation (resp.
endomorphism)  if $\mu(x\alpha y)=\mu(x)\alpha y+x\alpha\mu(y)$
(resp. $\mu(x\alpha y)=\mu(x)\alpha\mu(y)$) holds for all $x,y\in M$
and $\alpha\in \Gamma$.

Some conclusions concerning maps and identities in (semi-)prime
$\Gamma$-rings or $\Gamma_N$-rings  were obtained as a complement of
those in (semi-)prime rings. For example we can search some of them
in
\cite{Paul-Halder-2010-JPhySci-JordanLeft,Sapanci-Nakajima-1997-MJ-Complete-prime-gamma,
Soyturk-1994-TurkMJ-prime-gamma,Wang-1997-JNSM-prime-gamma,ZhangY-1991-JMRE-prime-gamma,
Paul-Halder-2009-JPhySci-JordanLeft,Ozturk-2002-TurkJM-der-prim-gamma,HuXQ-1998-thesis,Chakraborty-Paul-2010-BIMS}.

 In 2007, Asci and Ceran
\cite{Asci-Ceran-2007-IMForum-LeftPrimeGamma} discussed the left
derivation of a prime $\Gamma$-ring  and obtained some conclusions
on commutativity of a prime $\Gamma$-ring. For example it is proved
that for a prime $\Gamma$-ring $M$ with a nonzero ideal $U$ and a
nonzero right derivation $d$, if $charM\neq2,3$, $d^2(U)\subseteq
Z(M)$ and $d(U)\subseteq U$, then $M$ is commutative.

 In 2009, Paul
and Halder \cite{Paul-Halder-2009-JPhySci-JordanLeft} proved that a
left derivation of a semiprime $\Gamma$-ring $M$ must map $M$ into
its center under the assumption that $a\alpha b\beta c=a\beta
b\alpha c$ holds for all $a,b,c\in M$ and $\alpha,\beta\in \Gamma$.
In 2010, Paul and Halder \cite{Paul-Halder-2010-JPhySci-JordanLeft}
also considered this problem. In this paper, we get rid of the
assumptions in Paul and Halder's Theorems. Particularly we obtain
that a prime $\Gamma$-ring with a nonzero left or right derivation
must be commutative, which extends
 the conclusions appearing in
 \cite{Asci-Ceran-2007-IMForum-LeftPrimeGamma}.

 For a ring $R$ with $a,b\in R$, the symbol $[a,b]$ denotes
$ab-ba$. Similarly in a $\Gamma$-ring $M$ with $a,b\in M$ and
$\alpha\in \Gamma$, the symbol $[a,b]_{\alpha}$ denotes $a\alpha
b-b\alpha a$. For a ring $R$ with $a,b,c\in R$ the commutator
formulas $[a,bc]=[a,b]c+b[a,c]$ and $[ab,c]=a[b,c]+[a,c]b$ are very
useful in the course of dealing with identities in semiprime rings.
But in a $\Gamma$-ring $M$ with $a,b,c\in M$ and
$\alpha,\beta\in\Gamma$, the commutator formulas become $[a,b\alpha
c]_{\beta}=[a,b]_{\beta}\alpha c+b\alpha[a,c]_{\beta}+b\beta a\alpha
c-b\alpha a\beta c$ and $[a\alpha b,c]_{\beta}=[a,c]_{\beta}\alpha
b+a\alpha[b,c]_{\beta}+a\alpha c\beta b-a\beta c\alpha b$. The
assumption that  $a\alpha b\beta c=a\beta b\alpha c$ holds for all
$a,b,c\in M$ and $\alpha,\beta\in \Gamma$ can make the commutator
 formulas appear as  $[a,b\alpha c]_{\beta}=[a,b]_{\beta}\alpha
c+b\alpha[a,c]_{\beta}$ and $[a\alpha
b,c]_{\beta}=[a,c]_{\beta}\alpha b+a\alpha[b,c]_{\beta}$ which is
same as those in general rings. However in this paper, all results
being proved need not this assumption, which shows that some results
can be kept from semiprime rings to semiprime $\Gamma$-rings
although the basic commutator formulas have to be changed. For
convenience $a\alpha c\beta b-a\beta c\alpha b$ is usually denoted
by the symbol $a[\alpha,\beta]_cb$ although the form
$[\alpha,\beta]_c=\alpha c\beta- \beta c\alpha$ has not any sense.

A map $f$ from a ring $R$ into itself is called strong commutativity
preserving (scp) on a subset $S$ of $R$ if $[f(x),f(y)]=[x,y]$ holds
for all $x,y\in S$. The notion of a  strong commutativity preserving
map was first introduced by Bell and Mason
\cite{Bell-1994-MJOU-SCP}. Bell and Daif \cite{Bell-1994-CMB-SCP}
gave characterization of scp derivations and endomorphisms on
onesided ideals of semiprime rings. Bre\v{s}ar and Miers
\cite{Bresar-Robert-1994-CMB-scp} also described scp additive maps
on a semiprime ring.

Naturally a map $f$ from a $\Gamma$-ring $M$ into itself is called
strong commutativity preserving (scp) on a subset $S$ of $M$ if
$[f(x),f(y)]_{\alpha}=[x,y]_{\alpha}$ holds for all $x,y\in S$ and
$\alpha\in\Gamma$. In this paper, we also obtain that  a semiprime
$\Gamma$-ring with a strong commutativity preserving derivation on
itself must be commutative and that the strong commutativity
preserving endomorphism on a semiprime $\Gamma$-ring $M$ must be the
form $\sigma(x)=x+\zeta(x)$ where $\zeta$ is a map from $M$ into its
center $Z(M)$, which have ever been obtained  for semiprime rings by
Bell and Daif (\cite[Corollary 1]{Bell-1994-CMB-SCP} and the results
implied by \cite[Theorem 3]{Bell-1994-CMB-SCP}).

\section{Left derivations on semiprime $\Gamma$-rings}
The method employed  in this paper is straightforward computation.
Firstly we introduce some remarks.
\begin{rem} Let $M$ be a $\Gamma$-ring and $\delta:M\rightarrow M$ a left derivation.
Then  both $\delta([a,b]_{\alpha})=0$ and
 $[c,b]_{\beta}\alpha
\delta(a)=a\alpha c\beta\delta(b)-c\beta a\alpha\delta(b)$ hold for
\label{remark-leftD} all $a,b,c\in M$ and $\alpha,\beta\in \Gamma$.
\end{rem}

{\bf Proof.} Using the definition of a left derivation on a
$\Gamma$-ring and  $\delta((a\alpha b)\beta c)=\delta(a\alpha
(b\beta c))$  for all $a,b,c\in M$ and $\alpha,\beta\in \Gamma$, we
can obtain the conclusion through the straightforward computation.
\hfill $\Box$ \vskip0.3cm

The following simple observation is important for (semi-)prime
$\Gamma$-rings.

\begin{rem}
 Let $M$ be a $\Gamma$-ring with $c\in Z(M)$, $a_1,\ldots,a_n\in M$ and $\beta_1,\ldots,\beta_n\in
 \Gamma$.\nobreak Then  $$c\beta_1a_1\cdots\beta_na_n=a_1\beta_{\sigma(1)}\cdots
 a_i\beta_{\sigma(i)}c\beta_{{\sigma(i+1)}}a_{i+1}\cdots\beta_{\sigma(n)}a_n$$ for all $i\in\{1,\ldots,n\}$ and  $\sigma\in S_n$
  the symmetric group of degree $n$.\label{remark-center}
\end{rem}

{\bf Proof.} It is easy to see that for each $i\in \{1,2,\ldots,n\}$
 $$c\beta_1a_1\cdots\beta_na_n=a_1\beta_1\cdots
 a_i\beta_ic\beta_{i+1}a_{i+1}\cdots\beta_na_n$$
since $c\in Z(M)$. So it is sufficient to prove the equation
$$c\beta_1a_1\cdots\beta_na_n=c\beta_{\sigma(1)}a_1\cdots\beta_{\sigma(n)}a_n$$
for all $\sigma\in S_n$ in order to complete the proof.
 Now we prove this equation by
induction on $n$. For $n=1$ the equation is obvious. Set $n>1$.
Suppose that for every $k<n$ the equation holds. If $\sigma(n)=n$,
then there exists $\tau\in S_{n-1}$ such that $\tau(i)=\sigma(i)$
for all $1\leq i\leq n-1$. By the
 inductive assumption we get
 $$c\beta_1a_1\cdots\beta_{n-1}a_{n-1}=c\beta_{\tau(1)}a_1\cdots\beta_{\tau(n-1)}a_{n-1}=c\beta_{\sigma(1)}a_1\cdots\beta_{\sigma(n-1)}a_{n-1}.$$
Hence
$c\beta_1a_1\cdots\beta_na_n=c\beta_{\sigma(1)}a_1\cdots\beta_{\sigma(n)}a_n$.
If $\sigma(n)\neq n$, then by the
 inductive assumption  and $c\in Z(M)$ we
deduce that
 $$\begin{array}{rl}
 &c\beta_1a_1\cdots\beta_na_n\\
 =&(c\beta_{\sigma(n)}a_1\beta_{j_1}\cdots a_{n-2}\beta_{j_{n-2}}a_{n-1})\beta_na_n\\
 =&c\beta_{\sigma(n)}(a_1\beta_{j_1}\cdots a_{n-2}\beta_{j_{n-2}}a_{n-1}\beta_na_n)\\
 =&(a_1\beta_{j_1}\cdots a_{n-2}\beta_{j_{n-2}}a_{n-1}\beta_na_n)\beta_{\sigma(n)}c\\
 =&(c\beta_{j_1}a_1\cdots
 \beta_{j_{n-2}}a_{n-2}\beta_na_{n-1})\beta_{\sigma(n)}a_n\\
 =&c\beta_{\sigma(1)}a_1\cdots\beta_{\sigma(n-1)}a_{n-1}\beta_{\sigma(n)}a_n,
\end{array}$$
which completes the proof. \hfill $\Box$ \vskip0.3cm

 Remark \ref{remark-center} is important for prime or semiprime $\Gamma$-rings although it is
 easy to prove. For example if a semiprime $\Gamma$-ring $M$ has its center nonzero
 then every commutator formula has a neat form  with the help of any nonzero center
 element, i.e.,
 $d\alpha[a\beta b,c]_{\gamma}=d\alpha(a\beta[b,c]_{\gamma}+[a,c]_{\gamma}\beta
 b)$ holds for all $d\in Z(M)$, $a,b,c\in M$ and $\alpha,\beta,\gamma\in\Gamma$. Particularly for a prime $\Gamma$-ring
with its center nonzero every commutator formula has the same form
as one in a prime
 ring. That is $[a\beta b,c]_{\gamma}=a\beta[b,c]_{\gamma}+[a,c]_{\gamma}\beta
 b$ always holds for all $a,b,c\in M$ and $\beta,\gamma\in\Gamma$ in a prime $\Gamma$-ring with its center nonzero. But
 in general for most prime  or semiprime $\Gamma$-rings the center is zero. However
 Remark \ref{remark-center} is still useful for a semiprime $\Gamma$-ring when proving some results on commutativity even
 though the center is equal to zero. The  following characterization  for left
 derivations in prime or semiprime $\Gamma$-rings will make use of this
 observation.

\begin{thm}
A left  derivation  of a semiprime $\Gamma$-ring $M$ must map $M$
into its center. \label{theorem-main}
\end{thm}

{\bf Proof.} Let $\delta:M\rightarrow M$ be a left derivation. By
Remark \ref{remark-leftD} we have
\begin{align}
[c,b]_{\beta}\alpha\delta(a)=a\alpha c\beta\delta(b)-c\beta
a\alpha\delta(b),  ~~~~~~~~  a,b,c\in M, ~~~ \alpha, \beta\in
\Gamma.\label{equation-base}
\end{align}
 Putting
$b=[b,d]_{\gamma}$ in \eqref{equation-base} and applying Remark
\ref{remark-leftD} we obtain
\begin{align}
[c,[b,d]_{\gamma}]_{\beta}\alpha\delta(a)=0,~~~~~~~~ a,b,c,d\in M,
~~~ \alpha, \beta,\gamma\in \Gamma. \label{equation-second}
\end{align}
Then for all $a,b,c,d,a_1\in M$ and $\alpha,\beta,\gamma,\gamma_1\in
\Gamma$
\begin{align}
[a,[c,[b,d]_{\gamma}]_{\beta}]_{\alpha}\gamma_1\delta(a_1)=a\alpha[c,[b,d]_{\gamma}]_{\beta}\gamma_1\delta(a_1)=0.
\label{equation-third}
\end{align}
Putting $a=a\gamma_1a_1$ in \eqref{equation-second} and applying
\eqref{equation-third} we get that
$$\begin{array}{rcl}
0&=&[c,[b,d]_{\gamma}]_{\beta}\alpha
a_1\gamma_1\delta(a)+[c,[b,d]_{\gamma}]_{\beta}\alpha
a\gamma_1\delta(a_1)\\
&=&[c,[b,d]_{\gamma}]_{\beta}\alpha
a_1\gamma_1\delta(a)+[c,[b,d]_{\gamma}]_{\beta}\alpha
a\gamma_1\delta(a_1)+[a,[c,[b,d]_{\gamma}]_{\beta}]_{\alpha}\gamma_1\delta(a_1)\\
&=&[c,[b,d]_{\gamma}]_{\beta}\alpha
a_1\gamma_1\delta(a)+a\alpha[c,[b,d]_{\gamma}]_{\beta}\gamma_1\delta(a_1)\\
&=&[c,[b,d]_{\gamma}]_{\beta}\alpha a_1\gamma_1\delta(a)
\end{array}$$
holds for all $a,b,c,d,a_1\in M$ and
$\alpha,\beta,\gamma,\gamma_1\in \Gamma$. That is
$[c,[b,d]_{\gamma}]_{\beta}\Gamma M \Gamma\delta(a)=0$ holds for all
$a,b,c,d\in M$ and $\beta,\gamma\in \Gamma$. Hence
$[c,[b,\delta(a)]_{\gamma}]_{\beta}\Gamma M
\Gamma[c,[b,\delta(a)]_{\gamma}]_{\beta}=0$ holds for all $a,b,c\in
M$ and $\beta,\gamma\in \Gamma$. Then $[\delta(a),b]_{\gamma}\in
Z(M)$ for all $a,b\in M$ and $\gamma\in \Gamma$ since $M$ is
semiprime.

Put $a=d\gamma a$ in \eqref{equation-base} then for all $a,b,c,d\in
M$ and $\alpha,\beta,\gamma\in \Gamma$
\begin{align}
[c,b]_{\beta}\alpha d\gamma\delta(a)+[c,b]_{\beta}\alpha
a\gamma\delta(d)=d\gamma a\alpha c\beta\delta(b)-c\beta d\gamma
a\alpha\delta(b). \label{equation-cha1}
\end{align}
Multiply the two sides of \eqref{equation-base} with ``$d\gamma$"
from the left hand side then for all $a,b,c,d\in M$ and
$\alpha,\beta,\gamma\in \Gamma$
\begin{align}
d\gamma[c,b]_{\beta}\alpha\delta(a)=d\gamma a\alpha
c\beta\delta(b)-d\gamma c\beta a\alpha\delta(b).
\label{equation-cha2}
\end{align}
Compute \eqref{equation-cha1}$-$\eqref{equation-cha2} then for all
$a,b,c,d\in M$ and $\alpha,\beta,\gamma\in \Gamma$
\begin{align}
[c,b]_{\beta}\alpha d\gamma\delta(a)+[c,b]_{\beta}\alpha
a\gamma\delta(d)-d\gamma[c,b]_{\beta}\alpha\delta(a)=d\gamma c\beta
a\alpha\delta(b)-c\beta d\gamma a\alpha\delta(b).
\label{equation-charesult}
\end{align}
Setting $c=b$, $d=[\delta(b),b]_{\beta}\sigma\delta(b)$ in
\eqref{equation-charesult}, and then applying
$[\delta(b),b]_{\beta}\in Z(M)$ and Remark \ref{remark-center} we
have that
$$[\delta(b),b]_{\beta}\sigma[\delta(b),b]_{\beta}\gamma a\alpha\delta(b)=0$$
holds for all $a,b\in M$ and $\alpha,\beta,\gamma,\sigma\in \Gamma$.
That is $$[\delta(b),b]_{\beta}\sigma[\delta(b),b]_{\beta}\Gamma
M\Gamma[\delta(b),b]_{\beta}\sigma[\delta(b),b]_{\beta}=0$$ holds
for all $b\in M$ and $\beta,\sigma\in \Gamma$. Then
$[\delta(b),b]_{\beta}\sigma[\delta(b),b]_{\beta}=0$ holds for all
$b\in M$ and $\beta,\sigma\in \Gamma$. Hence
$[\delta(b),b]_{\beta}=0$ for all $b\in M$ and $\beta\in \Gamma$
since $[\delta(b),b]_{\beta}\in Z(M)$ the center of the semiprime
$\Gamma$ ring $M$. Setting $c=\delta(b)$,
$d=[d,\delta(b)]_{\beta}\sigma d$ in \eqref{equation-charesult}, and
then applying $[d,\delta(b)]_{\beta}\in Z(M)$ and Remark
\ref{remark-center} we deduce that
$$[d,\delta(b)]_{\beta}\sigma[d,\delta(b)]_{\beta}\gamma a\alpha\delta(b)=0$$
holds for all $a,b,d\in M$ and $\alpha,\beta,\gamma,\sigma\in
\Gamma$. Similar to proving that $[\delta(b),b]_{\beta}=0$ we also
obtain $[d,\delta(b)]_{\beta}=0$ for all $b,d\in M$ and
$\beta\in\Gamma$,
 which completes the proof.\hfill $\Box$
 \vskip0.3cm
Furthermore we will get the result for prime $\Gamma$-rings.

\begin{cor}
A prime $\Gamma$-ring with a  nonzero left derivation must be
commutative.
\end{cor}

{\bf Proof.} Let $\delta:M\rightarrow M$ be a nonzero left
derivation of a prime $\Gamma$-ring $M$. By Theorem
\ref{theorem-main}, Remark \ref{remark-leftD} and
\ref{remark-center}
\begin{align}
[c,b]_{\beta}\alpha\delta(a)=[a,c]_{\beta}\alpha\delta(b), ~~~~~~~~
a,b,c\in M, ~~~ \alpha, \beta\in \Gamma. \label{equation-prime-main}
\end{align}
Putting $c=a$ and applying $\delta(a)\in Z(M)$ we deduce that
$[a,b]_{\beta}\Gamma M\Gamma\delta(a)=0$ holds for all   $a,b\in M$
and $\beta\in \Gamma$. Hence for every $a\in M$ we deduce that
either $a\in Z(M)$ or $\delta(a)=0$. That is $M=\ker\delta\cup Z(M)$
is the union of its two subgroups. Thus  $M$ is commutative since
$\delta\neq0$.\hfill $\Box$

\section{Strong commutativity preserving maps on semiprime $\Gamma$-rings}

The following results (Theorem \ref{theorem-strong-der-semi} and
\ref{theorem-strong-endomo-semi}) on scp maps have been proved in
semiprime rings by Bell and Daif (see \cite{Bell-1994-CMB-SCP} for
reference in which more general situation were considered). Here we
will indicate that some results appearing in
\cite{Bell-1994-CMB-SCP} also hold in semiprime $\Gamma$-rings
although the commutator formulas have become complicated.

\begin{thm}
A semiprime $\Gamma$-ring with a strong commutativity preserving
derivation must be commutative. \label{theorem-strong-der-semi}
\end{thm}

{\bf Proof.} Suppose  that $M$ is a semiprime $\Gamma$-ring with a
strong commutativity preserving derivation $\delta$ on $M$. That is
$[\delta(x),\delta(y)]_{\alpha}=[x,y]_{\alpha}$ for all $x,y\in M$
and $\alpha\in\Gamma$. Then for all $x,y,z\in M$ and
$\alpha,\beta\in\Gamma$
$$[x\beta z,y]_{\alpha}=[\delta(x\beta z),\delta(y)]_{\alpha}=[\delta(x)\beta z,\delta(y)]_{\alpha}+[x\beta\delta(z),\delta(y)]_{\alpha}.$$
Moreover we get that
$$
\begin{array}{rl}
&x\beta[z,y]_{\alpha}+[x,y]_{\alpha}\beta z+x[\beta,\alpha]_yz\\
=&\delta(x)\beta[z,\delta(y)]_{\alpha}+[\delta(x),\delta(y)]_{\alpha}\beta z+\delta(x)[\beta,\alpha]_{\delta(y)}z+\\
 &x\beta[\delta(z),\delta(y)]_{\alpha}+[x,\delta(y)]_{\alpha}\beta\delta(z)+x[\beta,\alpha]_{\delta(y)}\delta(z)
\end{array}
$$
holds for all  $x,y,z\in M$ and $\alpha,\beta\in\Gamma$. That is for
all $x,y,z\in M$ and $\alpha,\beta\in\Gamma$
\begin{align}
x[\beta,\alpha]_yz=\delta(x)\beta[z,\delta(y)]_{\alpha}
+\delta(x)[\beta,\alpha]_{\delta(y)}z+[x,\delta(y)]_{\alpha}\beta\delta(z)+x[\beta,\alpha]_{\delta(y)}\delta(z).\label{equation-strong-derivation-base}
\end{align}
Putting $z=z\gamma t$ in \eqref{equation-strong-derivation-base} we
obtain that for all $x,y,z,t\in M$ and
$\alpha,\beta,\gamma\in\Gamma$
\begin{align}
\begin{array}{rcl}
x[\beta,\alpha]_yz\gamma
t&=&\delta(x)\beta\Big(z\gamma[t,\delta(y)]_{\alpha}+[z,\delta(y)]_{\alpha}\gamma
t+z[\gamma,\alpha]_{\delta(y)}t\Big)\\
&+&\delta(x)[\beta,\alpha]_{\delta(y)}z\gamma t+
[x,\delta(y)]_{\alpha}\beta\delta(z)\gamma t\\
&+&[x,\delta(y)]_{\alpha}\beta z\gamma\delta(t)
+x[\beta,\alpha]_{\delta(y)}\delta(z)\gamma
t+x[\beta,\alpha]_{\delta(y)}z\gamma\delta(t).
\end{array}\label{equation-strong-der-2}
\end{align}
Multiplying the two sides of \eqref{equation-strong-derivation-base}
by ``$\gamma t$" from the right hand side, and then comparing with
\eqref{equation-strong-der-2} we deduce that for all $x,y,z,t\in M$
and $\alpha,\beta,\gamma\in\Gamma$
\begin{align}
\delta(x)\beta\big(z\gamma[t,\delta(y)]_{\alpha}+z[\gamma,\alpha]_{\delta(y)}t\big)+
\big([x,\delta(y)]_{\alpha}\beta z
+x[\beta,\alpha]_{\delta(y)}z\big)\gamma\delta(t)=0.\label{equation-strong-der3}
\end{align}
Setting $t=\delta(y)$ and $\gamma=\alpha$ in
\eqref{equation-strong-der3} we get that for all  $x,y,z\in M$ and
$\alpha,\beta\in\Gamma$
\begin{align}
(x\beta\delta(y)\alpha z-\delta(y)\alpha x\beta
z)\alpha\delta^2(y)=0. \label{equation-strong-der-last1}
\end{align}
Putting $\alpha=\alpha+\gamma$ into
\eqref{equation-strong-der-last1} and applying
\eqref{equation-strong-der-last1} we deduce that for all $x,y,z\in
M$ and $\alpha,\beta,\gamma\in\Gamma$
$$(x\beta\delta(y)\alpha z-\delta(y)\alpha x\beta
z)\gamma\delta^2(y)=-(x\beta\delta(y)\gamma z-\delta(y)\gamma x\beta
z)\alpha\delta^2(y).$$ Then by \eqref{equation-strong-der-last1} for
all $x,y,z\in M$ and $\alpha,\beta,\gamma\in\Gamma$
$$
\begin{array}{ll}
&(x\beta\delta(y)\alpha z-\delta(y)\alpha x\beta
z)\gamma\delta^2(y)\Gamma M\Gamma(x\beta\delta(y)\alpha
z-\delta(y)\alpha x\beta z)\gamma\delta^2(y)\\
=-&(x\beta\delta(y)\alpha z-\delta(y)\alpha x\beta
z)\gamma\delta^2(y)\Gamma M\Gamma(x\beta\delta(y)\gamma
z-\delta(y)\gamma x\beta z)\alpha\delta^2(y)\\
=0,&
\end{array}
$$
which implies $(x\beta\delta(y)\alpha z-\delta(y)\alpha x\beta
z)\gamma\delta^2(y)=0$ for all $x,y,z\in M$ and
$\alpha,\beta,\gamma\in\Gamma$ since $M$ is semiprime. Set
$\beta=\alpha$ and $x=\delta(x)$ in $(x\beta\delta(y)\alpha
z-\delta(y)\alpha x\beta z)\gamma\delta^2(y)=0$. Then for all
$x,y,z\in M$ and $\alpha\in\Gamma$
$$
\begin{array}{rl}
&[x,y]_{\alpha}\alpha z\Gamma
M\Gamma[x,y]_{\alpha}\alpha z\\
=&[\delta(x),\delta(y)]_{\alpha}\alpha z\Gamma
M\Gamma[\delta(x),\delta(y)]_{\alpha}\alpha z\\
=&[\delta(x),\delta(y)]_{\alpha}\alpha z\Gamma
M\Gamma[\delta^2(x),\delta^2(y)]_{\alpha}\alpha z\\
=&0,
\end{array}
$$
which shows $[x,y]_{\alpha}\alpha z=0$ for all $x,y,z\in M$ and
$\alpha\in\Gamma$. Then for all $t,x,y,z\in M$ and $\alpha,\gamma\in
\Gamma$
\begin{align}
0=[t\gamma x,y]_{\alpha}\alpha z=t\gamma[x,y]_{\alpha}\alpha
z+(t\gamma y\alpha x-y\alpha t\gamma x)\alpha z=(t\gamma y\alpha
x-y\alpha t\gamma x)\alpha z. \label{equation-strong-der-last2}
\end{align}
Putting $\alpha=\alpha+\beta$ into \eqref{equation-strong-der-last2}
and applying \eqref{equation-strong-der-last2} we deduce that for
all $t,x,y,z\in M$ and $\alpha,\beta,\gamma\in\Gamma$
$$(t\gamma y\alpha x-y\alpha t\gamma x)\beta z=-(t\gamma y\beta x-y\beta t\gamma x)\alpha z.$$ Then by
\eqref{equation-strong-der-last2} for all $t,x,y,z\in M$ and
$\alpha,\beta,\gamma\in\Gamma$
$$
\begin{array}{ll}
&(t\gamma y\alpha x-y\alpha t\gamma x)\beta z\Gamma
M\Gamma(t\gamma y\alpha x-y\alpha t\gamma x)\beta z\\
=-&(t\gamma y\alpha x-y\alpha t\gamma x)\beta z\Gamma
M\Gamma(t\gamma y\beta x-y\beta t\gamma x)\alpha z\\
=0,&
\end{array}
$$
which implies $(t\gamma y\alpha x-y\alpha t\gamma x)\beta z=0$ for
all $t,x,y,z\in M$ and $\alpha,\beta,\gamma\in\Gamma$. Moreover
$t\gamma y\alpha x-y\alpha t\gamma x=0$ for all $t,x,y\in M$ and
$\alpha,\gamma\in\Gamma$. Hence for all  $t,x,y\in M$ and
$\alpha,\gamma\in\Gamma$
$$[t\gamma x, y]_{\alpha}=t\gamma x\alpha y-y\alpha t\gamma x=t\gamma x\alpha y-t\gamma y\alpha x=t\gamma[x,y]_{\alpha}.$$
Then for all $x,y,z\in M$ and $\alpha,\beta\in\Gamma$
$$
\begin{array}{rl}
&x\beta[z,y]_{\alpha}=[x\beta z,y]_{\alpha}=[\delta(x\beta
z),\delta(y)]_{\alpha}\\
=&[\delta(x)\beta
z,\delta(y)]_{\alpha}+[x\beta\delta(z),\delta(y)]_{\alpha}\\
=&\delta(x)\beta[z,\delta(y)]_{\alpha}+x\beta[\delta(z),\delta(y)]_{\alpha}.
\end{array}
$$
So $\delta(x)\beta[z,\delta(y)]_{\alpha}=0$ for all $x,y,z\in M$ and
$\alpha,\beta\in\Gamma$. Hence for all $t,x,y,z\in M$ and
$\alpha,\beta,\gamma\in\Gamma$
 $$0=\delta(x)\beta[t\gamma z,\delta(y)]_{\alpha}=\delta(x)\beta t\gamma[z,\delta(y)]_{\alpha},$$
which implies $[\delta(x),\delta(y)]_{\alpha}\Gamma
M\Gamma[\delta(x),\delta(y)]_{\alpha}=0$ for all $x,y\in M$ and
$\alpha\in\Gamma$. Thus
$[x,y]_{\alpha}=[\delta(x),\delta(y)]_{\alpha}=0$ for all $x,y\in M$
and $\alpha\in\Gamma$ completes the proof. \hfill $\Box$

\vskip0.3cm It is implied by \cite[Theorem 3]{Bell-1994-CMB-SCP}
 that for  a semiprime ring  $R$ with an
endomorphism $T$, then $T$ is scp on $R$ if and only if
$T(x)=x+\zeta(x)$ for all $x\in R$ where $\zeta$ is a map from $R$
into its center. We will show this also holds for semiprime
$\Gamma$-rings.

\begin{thm}
Let $M$ be a semiprime $\Gamma$-ring with an endomorphism $\sigma$.
Then $\sigma$ is strong commutativity preserving on $M$ if and only
if there exists a map $\zeta:M\rightarrow Z(M)$ such that
$\sigma(x)=x+\zeta(x)$ for all $x\in M$.
\label{theorem-strong-endomo-semi}
\end{thm}

{\bf Proof.} We will only consider the necessity since the
sufficiency is obvious. From $[\sigma(x\alpha z),
\sigma(x)]_{\alpha}=[x\alpha z,x]_{\alpha}$ for all $x,z\in M$ and
$\alpha\in\Gamma$ we obtain $(\sigma(x)-x)\alpha[z,x]_{\alpha}=0$.
Then for all $x,y,z\in M$ and $\alpha,\beta\in\Gamma$
\begin{align}
\begin{array}{rcl}
0&=&(\sigma(x)-x)\alpha[y\beta
z,x]_{\alpha}\\
&=&(\sigma(x)-x)\alpha[y,x]_{\alpha}\beta
z+(\sigma(x)-x)\alpha(y\beta z\alpha x-y\alpha x\beta z)\\
&=&(\sigma(x)-x)\alpha(y\beta z\alpha x-y\alpha x\beta z).
\end{array}\label{equation-strong-auto1}
\end{align}
Linearizing $\alpha$ in \eqref{equation-strong-auto1} we have that
for all $x,y,z\in M$ and $\alpha,\beta,\gamma\in\Gamma$
$$(\sigma(x)-x)\alpha(y\beta z\gamma x-y\gamma x\beta z)=-(\sigma(x)-x)\gamma(y\beta z\alpha x-y\alpha x\beta z).$$
Then by \eqref{equation-strong-auto1} for all $x,y,z\in M$ and
$\alpha,\beta,\gamma\in\Gamma$
$$
\begin{array}{rcl}
&&(\sigma(x)-x)\alpha(y\beta z\gamma x-y\gamma x\beta z)\Gamma
M\Gamma(\sigma(x)-x)\alpha(y\beta z\gamma x-y\gamma x\beta z)\\
&=&-(\sigma(x)-x)\alpha(y\beta z\gamma x-y\gamma x\beta z)\Gamma
M\Gamma(\sigma(x)-x)\gamma(y\beta z\alpha x-y\alpha x\beta z)\\
&=&0
\end{array}
$$
which implies $(\sigma(x)-x)\alpha(y\beta z\gamma x-y\gamma x\beta
z)=0$ for all $x,y,z\in M$ and $\alpha,\beta,\gamma\in\Gamma$.
Linearizing $x$ in $(\sigma(x)-x)\alpha(y\beta z\gamma x-y\gamma
x\beta z)=0$ we get for all  $x,y,z,t\in M$ and
$\alpha,\beta,\gamma\in\Gamma$.
$$(\sigma(x)-x)\alpha(y\beta z\gamma t-y\gamma
t\beta z)=-(\sigma(t)-t)\alpha(y\beta z\gamma x-y\gamma x\beta z).$$
Then for all  $x,y,z,t\in M$ and $\alpha,\beta,\gamma\in\Gamma$
$$
\begin{array}{rcl}
&&(\sigma(x)-x)\alpha(y\beta z\gamma t-y\gamma t\beta z)\Gamma
M\Gamma(\sigma(x)-x)\alpha(y\beta z\gamma t-y\gamma
t\beta z)\\
&=&-(\sigma(x)-x)\alpha(y\beta z\gamma t-y\gamma t\beta z)\Gamma
M\Gamma(\sigma(t)-t)\alpha(y\beta z\gamma x-y\gamma x\beta z)\\
&=&0,
\end{array}
$$
which shows $(\sigma(x)-x)\alpha(y\beta z\gamma t-y\gamma t\beta
z)=0$ for all $x,y,z,t\in M$ and $\alpha,\beta,\gamma\in\Gamma$.
Then
 for all $x,y,z\in M$ and $\alpha,\beta,\gamma\in\Gamma$ we obtain both
$$((\sigma(y)-y)\beta z\gamma x-(\sigma(y)-y)\gamma x\beta
z)\Gamma M\Gamma((\sigma(y)-y)\beta z\gamma x-(\sigma(y)-y)\gamma
x\beta z)=0$$ and
$$(y\beta (\sigma(x)-x)\gamma z-y\gamma z\beta
(\sigma(x)-x))\Gamma M\Gamma(y\beta (\sigma(x)-x)\gamma z-y\gamma
z\beta (\sigma(x)-x))=0.$$ So both $(\sigma(y)-y)\beta z\gamma
x-(\sigma(y)-y)\gamma x\beta z=0$ and $y\beta (\sigma(x)-x)\gamma
z-y\gamma z\beta (\sigma(x)-x)=0$ hold for all $x,y,z\in M$ and
$\beta,\gamma\in\Gamma$. Hence for all $x,y,z\in M$ and
$\alpha,\beta\in\Gamma$ we get
$$[(\sigma(x)-x)\alpha y,z]_{\beta}=[\sigma(x)-x,z]_{\beta}\alpha
y,
~~[y\alpha(\sigma(x)-x),z]_{\beta}=[y,z]_{\beta}\alpha(\sigma(x)-x)$$
and $[y\alpha z,\sigma(x)-x]_{\beta}=[y,\sigma(x)-x]_{\beta}\alpha
z$. Then for all $x,y,z\in M$ and $\alpha,\beta\in\Gamma$
$$
\begin{array}{rcl}
0&=&[\sigma(x)\alpha\sigma(y),\sigma(z)]_{\beta}-[x\alpha
y,z]_{\beta}\\
&=&\big([\sigma(x)\alpha\sigma(y),\sigma(z)]_{\beta}-[\sigma(x)\alpha
y,\sigma(z)]_{\beta}\big)+\big([\sigma(x)\alpha
y,\sigma(z)]_{\beta}-[x\alpha
y,\sigma(z)]_{\beta}\big)\\
&&+\big([x\alpha y,\sigma(z)]_{\beta}-[x\alpha y,z]_{\beta}\big)\\
&=&[\sigma(x)\alpha(\sigma(y)-y),\sigma(z)]_{\beta}+[(\sigma(x)-x)\alpha
y,\sigma(z)]_{\beta}+[x\alpha y,\sigma(z)-z]_{\beta}\\
&=&[\sigma(x),\sigma(z)]_{\beta}\alpha(\sigma(y)-y)+[\sigma(x)-x,\sigma(z)]_{\beta}\alpha
y+[x,\sigma(z)-z]_{\beta}\alpha y\\
&=&[x,z]_{\beta}\alpha(\sigma(y)-y)+[x,z]_{\beta}\alpha
y-[x,\sigma(z)]_{\beta}\alpha y+[x,\sigma(z)-z]_{\beta}\alpha y\\
&=&[x,z]_{\beta}\alpha(\sigma(y)-y).
\end{array}
$$
Hence for all $x,y,z\in M$ and $\alpha,\beta,\gamma\in\Gamma$
$$0=[(\sigma(y)-y)\gamma x,z]_{\beta}\alpha(\sigma(y)-y)=[\sigma(y)-y,z]_{\beta}\gamma x\alpha(\sigma(y)-y).$$
Thus $[\sigma(y)-y,z]_{\beta}\Gamma
M\Gamma[\sigma(y)-y,z]_{\beta}=0$ holds for all $y,z\in M$ and
$\beta\in\Gamma$. Then $[\sigma(y)-y,z]_{\beta}=0$ for all $y,z\in
M$ and $\beta\in\Gamma$ completes the proof. \hfill $\Box$

\vskip0.3cm For prime $\Gamma$-rings we get a further  result.

\begin{cor}
In a noncommutative  prime $\Gamma$-ring $M$ the identity map is the
unique strong commutativity preserving endomorphism on $M$.
\end{cor}

{\bf Proof.} Let $\sigma:M\rightarrow M$ be a strong commutativity
preserving endomorphism on $M$. Then by Theorem
\ref{theorem-strong-endomo-semi} there exists a map
$\zeta:M\rightarrow Z(M)$ such that $\sigma(x)=x+\zeta(x)$ for all
$x\in M$. For all $x,y\in M$ and $\alpha\in\Gamma$
\begin{align}
x\alpha y+\zeta(x\alpha y)=\sigma(x\alpha
y)=\sigma(x)\alpha\sigma(y)=x\alpha y+\zeta(x)\alpha
y+\zeta(y)\alpha
x+\zeta(x)\alpha\zeta(y).\label{equation-strong-endo-prime-main}
\end{align}
Then for $x\in M$ such that $\zeta(x)\neq 0$ we get
$\zeta(x)\alpha[x,y]_{\beta}=0$ for all $y\in M$ and
$\alpha,\beta\in\Gamma$ from \eqref{equation-strong-endo-prime-main}
and Remark \ref{remark-center}. So $[x,y]_{\beta}=0$ for all $y\in
M$ and $\beta\in\Gamma$ since $M$ is prime and $0\neq\zeta(x)\in
Z(M)$. That is for every $x\in M$ once $\zeta(x)\neq0$ we always
have $x\in Z(M)$. Now we assume that there exists $x_0\in M$ such
that $\zeta(x_0)\neq0$ and proceed to obtain a contradiction so that
the proof could be completed. We may choose an element $y_0\in
M\backslash Z(M)$ since $M$ is noncommutative. Then $\zeta(y_0)=0$,
i.e. $\sigma(y_0)=y_0$. Putting $x=x_0$ and $y=y_0$ in
\eqref{equation-strong-endo-prime-main} we obtain $\zeta(x_0)\alpha
y_0\in Z(M)$. So $\zeta(x_0)\alpha[y_0,z]_{\beta}=0$ holds for all
$z\in M$ and $\alpha,\beta\in\Gamma$ which means $y_0\in Z(M)$ a
contradiction. \hfill $\Box$

\vskip0.3cm Next we will give an example showing that there exists a
non-identity strong commutativity preserving endomorphism on certain
noncommutative semiprime $\Gamma$-ring.

{\bf Example}  Let $R=\mathbb{M}_2(\mathbb{C})\times \mathbb{C}$
where $\mathbb{C}$ is the field of complex numbers. Then $R$ is a
noncommutative semiprime ring and $\sigma:R\rightarrow R$ such that
$\sigma(A,a)=(A,\overline{a})$ for all $(A,a)\in R$ is a
non-identity strong commutativity preserving automorphism on
$R$.\hfill $\Box$

 {\bf Acknowledgement} We would like to thank Professor Paul for
 sending the paper required by us during  preparing  this
manuscript.

\bibliographystyle{amsplain}

\begin{thebibliography}{150}
{\baselineskip 0.668cm

\bibitem{Asci-Ceran-2007-IMForum-LeftPrimeGamma}Asci, M. and Ceran, S., {\it  The commutativity in prime $\Gamma$-rings with left derivation},
    Int. Math. Forum, {\bf 2}(2007), 103--108.

\bibitem{Barnes-1966-PJM-GammaRing}Barnes, W.E., {\it  On the $\Gamma$-rings of Nobusawa},
    Pacific J. Math., {\bf 18}(1966), 411--422.

\bibitem{Bell-1994-CMB-SCP}Bell, H.E. and Daif, M.N., {\it  On commutativity and strong commutativity preserving maps},
    Canad. Math. Bull., {\bf 37}(1994), 443--447.

\bibitem{Bell-1994-MJOU-SCP}Bell, H.E. and Mason, G., {\it  On derivations in near rings and rings},
    Math. J. Okayama Univ., {\bf 34}(1992), 135--144.


\bibitem{Bresar-Robert-1994-CMB-scp}Bre\v{s}ar, M. and Miers, C.R., {\it  Strong commutativity preserving maps of semiprime
rings},
    Canad. Math. Bull., {\bf 37}(1994), 457--460.

\bibitem{Bresar-Vukman-1990-PAMS-left}Bre\v{s}ar, M. and Vukman, J., {\it  On left derivations and related mappings},
    Proc. Amer. Math. Soc., {\bf 110}(1990), 7--16.

\bibitem{Chakraborty-Paul-2010-BIMS}Chakraborty, S. and Paul, A.C., {\it  On Jordan generalized $k$-derivations of semiprime $\Gamma_N$-rings},
    Bull. Iranian Math. Soc., {\bf 36}(2010), 41--53.


\bibitem{HuXQ-1998-thesis}Hu, X., {\it Derivations in prime gamma rings},
   Thesis (Ph.D.)--North Carolina State University, 1998.



\bibitem{Kyuno-1978-PJM-PrimeGammaRing}Kyuno, S., {\it  On prime $\Gamma$-rings},
    Pacific J. Math., {\bf 75}(1978), 185--190.

\bibitem{Luh-1969-MichiganMJ-Simple-Gamma}Luh, J., {\it  On the theory of simple $\Gamma$-rings},
    Michigan Math. J., {\bf 16}(1969), 65--75.


\bibitem{Nobusawa-1964-OsakaJM-GammaNRing}Nobusawa, N., {\it  On a generalization of the ring theory},
   Osaka J. Math., {\bf 1}(1964), 81--89.


\bibitem{Ozturk-2002-TurkJM-der-prim-gamma}Ozturk, M.A., Jun, Y.B. and Kim, K.H., {\it  On derivations of prime gamma rings},
   Turkish J. Math., {\bf 26}(2002), 317--327.




\bibitem{Paul-Halder-2009-JPhySci-JordanLeft}Paul, A.C. and Halder, A.K., {\it  Jordan left derivations of two torsion free $\Gamma M$-Modules},
  J. Phys. Sci., {\bf 13}(2009), 13--19.

\bibitem{Paul-Halder-2010-JPhySci-JordanLeft}Paul, A.C. and Halder, A.K., {\it On left derivations of $\Gamma$-rings},
  Bull. Pure Appl. Math., {\bf 4}(2010), 320--328.


\bibitem{Sapanci-Nakajima-1997-MJ-Complete-prime-gamma}Sapanci, M. and Nakajima, A., {\it Jordan derivations on completely
prime gamma rings},
  Math. Japon., {\bf 46}(1997), 47--51.

\bibitem{Soyturk-1994-TurkMJ-prime-gamma}Soyturk, M., {\it The commutativity in prime gamma rings with
derivation},
  Turkish J. Math., {\bf 18}(1994), 149--155.


\bibitem{Wang-1997-JNSM-prime-gamma}Wang, D., {\it On derivations and commutativity for prime
$\Gamma$-rings},
  J. Natur. Sci. Math., {\bf 37}(1997), 99--112.

\bibitem{ZhangY-1991-JMRE-prime-gamma}Zhang, Y., {\it Derivations in a prime $\Gamma$-ring},
 J. Math. Res. Exposition, {\bf 11}(1991), 565--568.


 }
\end{thebibliography}

\end{document}